\documentclass[11pt]{article}
\usepackage{graphicx}
 \usepackage{latexsym}
\usepackage{amssymb}

\begin{document}

\fontsize{11}{14.5pt}\selectfont

\begin{center}

{\small Technical Report No.\ 0411,
 Department of Statistics, University of Toronto}

\vspace*{0.9in}

{\Large\bf Taking Bigger Metropolis Steps by Dragging Fast Variables}\\[16pt]

{\large Radford M. Neal}\\[3pt]
 Department of Statistics and Department of Computer Science \\
 University of Toronto, Toronto, Ontario, Canada \\
 \texttt{http://www.cs.utoronto.ca/$\sim$radford/} \\
 \texttt{radford@stat.utoronto.ca}\\[10pt]

 26 October 2004 (some typos fixed 6 February 2005)

\end{center}

\vspace*{8pt}

\noindent \textbf{Abstract.}  I show how Markov chain sampling with
the Metropolis-Hastings algorithm can be modified so as to take bigger
steps when the distribution being sampled from has the characteristic
that its density can be quickly recomputed for a new point if this
point differs from a previous point only with respect to a subset of
``fast'' variables.  I show empirically that when using this method,
the efficiency of sampling for the remaining ``slow'' variables can
approach what would be possible using Metropolis updates based on the
marginal distribution for the slow variables.

\section{\hspace*{-7pt}Introduction}\label{sec-intro}\vspace*{-10pt}

Suppose we wish to sample from a distribution $\pi(x,y) \propto
\exp(-E(x,y))$, where $E$ is a given ``energy'' function, by
simulating a Markov chain with $\pi$ as its equilibrium distribution.
Let's suppose that $x$ is a ``slow'' variable and $y$ is a ``fast''
variable, so that once $E(x,y)$ has been computed (and intermediate
quantities cached), we can compute $E(x,y')$ much faster than we can
compute $E(x',y')$ for some $x'$ for which we haven't previously 
calculated $E$.

I was led to consider this problem because it arises with Bayesian
models that attempt to infer cosmological parameters from data on the
cosmic microwave background radiation (Lewis and Bridle 2002), for
which recomputation after changing only fast variables can be around a
thousand times faster than recomputation after changing a slow
variable.  Similarly large differences between fast and slow variables
can arise with Gaussian process classification models (Neal 1999), 
in which updating the latent variables is fast, while updating the
parameters of the covariance matrix is slow, since the new covariance
matrix must then be inverted.  Computationally equivalent problems
also arise in geostatistics (Diggle, Tawn, and Moyeed 1998), and for
what are called ``generalized linear mixed effects models''.
Many other statistical problems also have some variables that are
faster than others, though not always by such a large factor. 

Ideally, we would like to do Metropolis-Hastings updates for $x$ (Metropolis,
\textit{et al} 1953; Hastings 1970), using some
proposal distribution, $S(x^*|x)$, and accepting or rejecting $x^*$
based on its marginal distribution, $\pi(x)$.  The acceptance
probability for such a proposal would be
\begin{eqnarray}
  a(x,x^*) & = & \min\left[\, 1,\ 
  { S(x|x^*)\, \pi(x^*) \over S(x^*|x)\, \pi(x) } \right] 
\end{eqnarray}
Suppose, however, that we can't feasibly compute the
marginal distribution, $\pi(x)$, so that this approach is not possible.  
Instead we will have to use a Metropolis-Hastings algorithm that operates
on the joint distribution for $x$ and $y$.
If we could sample directly from the
conditional distribution for $y$, $\pi(y|x)$, we could generate $x^*$
from $S(x^*|x)$ and then
$y^*$ from $\pi(y^*|x^*)$, and the resulting acceptance probability
for $(x^*,y^*)$ would be the same (due to cancellation) as that above
using the marginal distribution for $x$.  However, let's assume that sampling
from $\pi(y|x)$ is also infeasible.  We might hope to approximate $\pi(y|x)$ by
some transition distribution $T(y^*|y;x)$ that we can sample from.  To use
this approximation in a Metropolis-Hasting proposal, however, we would
need to be able to compute the probability of proposing $y^*$, which
will likely be impossible if we have to resort to iterative methods
(eg, Markov chain simulation) in order to obtain a good approximation.

This paper describes a way in which these problems can be bypassed
when recomputing $E(x,y)$ after changing only the ``fast' variable $y$
is much quicker than recomputing $E(x,y)$ after changing $x$. 
In this method, changes to $x$ are made in conjunction with
changes to $y$ that are found by ``dragging'' $y$ with the help
of intermediate transitions that involve only fast re-computations
of $E$.  In the limit as the number of such intermediate transitions
increases, I show empirically (but haven't proved) that the method is equivalent
to using the marginal distribution of $x$.  Since the intermediate
transitions involve only fast computations, we hope to be able to do
quite a few intermediate transitions, and get close to the effect of
using the marginal probabilities for $x$.

The method can be seen as a generalization of ``tempered transitions''
(Neal 1996), and could be expressed in greater generality than I have
done here, where I concentrate on the context with fast and slow variables.
To begin, I'll describe the method when there is only one intermediate
transition, since this is easier to work with, but I expect that one
would use many intermediate transitions in practice, as described
later.

\section{\hspace*{-7pt}The method with one intermediate 
                       transition}\label{sec-one}\vspace{-10pt}

If the current state is $(x,y)$, we start by generating a proposed 
new value $x^*$ according to the probabilities $S(x^*|x)$.  We then
define a distribution, $\rho$, over values for $y$ that is intermediate 
between $\pi(y|x)$ and $\pi(y|x^*)$, as follows:
\begin{eqnarray}
    \rho(y;x,x^*) & \propto & \exp(-(E(x,y)+E(x^*,y))/2)
\label{eq-rho}
\end{eqnarray}
Here, the dependence of $\rho$ on $x$ and $x^*$ has been made explicit,
but note that this is a distribution over $y$ only, not $x$ and $y$
jointly.
We choose some transition probabilities, $T$, for updating $y$ so as to leave
$\rho$ invariant.  These probabilities must of course depend on $x$ and $x^*$.
We write them as $T(y'|y;x,x^*)$.  We require that they satisfy detailed
balance, so that for all $x$, $x^*$, $y$, and $y'$,
\begin{eqnarray}
  \rho(y;x,x^*)\,T(y'|y;x,x^*) & = & \rho(y';x,x^*)\,T(y|y';x,x^*)
\label{eq-db}
\end{eqnarray}
We also require that $T$ depend symmetrically on the two $x$ values ---
for all $x$, $x^*$, $y$ and $y'$:
\begin{eqnarray}
   T(y'|y;x,x^*) & = & T(y'|y;x^*,x)
\label{eq-Tsym}
\end{eqnarray}
$T$ might, for example, be a Metropolis-Hastings update, or a series
of such updates.  We apply
this transition once, to sample a value $y^*$ from $T(y^*|y;x,x^*)$.
We then accept $(x^*,y^*)$ as the next state with 
probability $a(x,y,x^*,y^*)$, defined as follows:
\begin{eqnarray}
  a(x,y,x^*,y^*) & = & \min\left[\, 1,\ 
  { S(x|x^*)\, \pi(x^*,y^*)\, \rho(y;x,x^*) 
    \over S(x^*|x)\, \pi(x,y)\, \rho(y^*;x,x^*)}\,
  \right] \\[4pt]
  & = & \min\left[\, 1,\ 
  { S(x|x^*)\, \exp(-E(x^*,y^*))\, \exp(-(E(x,y)+E(x^*,y))/2)
    \over S(x^*|x)\, \exp(-E(x,y))\, \exp(-(E(x,y^*)+E(x^*,y^*))/2)}\,
  \right]\ \ \ \ \\[4pt]
  & = & \min\left[\, 1,\ 
  { S(x|x^*) \over S(x^*|x)}\, 
  \exp\left({E(x,y)+E(x,y^*)\over2}\ -\ {E(x^*,y)+E(x^*,y^*) \over2} \right)
  \right]\ \ \ \ 
\end{eqnarray}
If we don't accept, the next state is the current state, $(x,y)$.

Although this expression for $a(x,y,x^*,y^*)$ has four occurrences of
$E(\cdot,\cdot)$, only two slow evaluations are needed.  In fact, only
one slow evaluation is needed if we assume that an evaluation was
done previously for the current state, when it was proposed.
Note also that we would often choose a symmetric proposal distribution 
for $x$, so that $S(x^*|x)/S(x|x^*)\,=\,1$.

To show that this is a valid update, I will prove that it satisfies 
detailed balance.  The probability of moving from $(x,y)$ to 
a different state $(x^*,y^*)$ when in equilibrium is\vspace*{-4pt}
\begin{eqnarray}
\lefteqn{
  \pi(x,y)\,S(x^*|x)\,T(y^*|y;x,x^*)\,
  \min\left[\, 1,\ 
  { S(x|x^*)\, \pi(x^*,y^*)\, \rho(y;x,x^*) 
    \over S(x^*|x)\, \pi(x,y)\, \rho(y^*;x,x^*)}\,
  \right]} \ \nonumber\\[4pt]
  & = &
  \min\left[\,S(x^*|x)\,\pi(x,y)\,T(y^*|y;x,x^*),\
  {S(x|x^*)\, \pi(x^*,y^*)\, \rho(y;x,x^*)\,T(y^*|y;x,x^*)
   \over \rho(y^*;x,x^*) }\,
  \right]\ \ \ \ \ \\[4pt]
  & = &
  \min\Big[\,S(x^*|x)\,\pi(x,y)\,T(y^*|y;x,x^*),\
   S(x|x^*)\, \pi(x^*,y^*)\, T(y|y^*;x,x^*)\, \Big]
   \ \ \ \ \ \\[4pt]
  & = &  
  \min\Big[\,S(x^*|x)\,\pi(x,y)\,T(y^*|y;x,x^*),\
   S(x|x^*)\, \pi(x^*,y^*)\, T(y|y^*;x^*,x)\, \Big]
\end{eqnarray}
Here, the detailed balance condition (\ref{eq-db}) and symmetry 
condition (\ref{eq-Tsym}) have been used.  Examination of the
above shows that swapping $(x,y)$ and $(x^*,y^*)$ leaves it
unchanged, showing the detailed balance holds.  

I would expect this method to work better than the simple method of
just proposing to change from $x$ to $x^*$ while keeping $y$ unchanged.
The latter method will work well only if the old $y$ is often suitable
for the new $x^*$ --- ie, if the old $y$ is typical of $\pi(y|x^*)$.
This will often be true only if the change from $x$ to $x^*$ is small.
The new method changes $y$ to a $y^*$ that is drawn approximately (if
$T$ works well) from a distribution that is halfway between $\pi(y|x)$
and $\pi(y|x^*)$.  Such a $y^*$ should have a better chance of being
suitable for $x^*$, allowing the change from $x$ to $x^*$ to be
greater while still maintaining a good acceptance probability.
If we propose an $x^*$ that is a really big change from $x$, however,
even a $y^*$ that comes from a distribution halfway to $\pi(y|x^*)$
may not be good enough.  

\section{\hspace*{-7pt}The method with many intermediate 
                          transitions}\label{sec-many}\vspace{-10pt}

We can try to take bigger steps in $x$ by ``dragging'' $y$
through a series of intermediate distributions interpolating between
$\pi(y|x)$ and $\pi(y|x^*)$.  Given some integer $n>1$, we define the 
following distributions, for $i=0,\ldots,n$:
\begin{eqnarray}
  \rho_i(y;x,x^*) & \propto & 
    \exp(\,-\,((1\!-\!i/n)E(x,y)\,+\,(i/n)E(x^*,y)))
\label{eq-rhoi}
\end{eqnarray}
Notice that $\rho_0(y;x,x^*)=\pi(y|x)$ and $\rho_n(y;x,x^*)=\pi(y|x^*)$.
When $n=2$, $\rho_1$ is the same as the $\rho$ defined above in (\ref{eq-rho}).
Finally, note that $\rho_i(y;x,x^*) = \rho_{n-i}(y;x^*,x)$.

For each $\rho_i$, we need to choose transition probabilities, $T_i$,
which may depend on $x$ and $x^*$.  We require that they satisfy detailed
balance, so that for all $x$, $x^*$, $y$, and $y'$,
\begin{eqnarray}
  \rho_i(y;x,x^*)\,T_i(y'|y;x,x^*) & = & \rho_i(y';x,x^*)\,T_i(y|y';x,x^*)
\label{eq-db2}
\end{eqnarray}
We also require of each opposite pair of transitions, $T_i$ and $T_{n-i}$,
that for all $x$, $x^*$, $y$ and $y'$,
\begin{eqnarray}
   T_i(y'|y;x,x^*) & = & T_{n-i}(y'|y;x^*,x)
\label{eq-Tsym2}
\end{eqnarray}
These conditions will be satisfied if the $T_i$ are standard 
Metropolis updates with respect to the $\rho_i$, with $T_i$ and $T_{n-i}$
using the same proposal distribution.

The update procedure using $n-1$ intermediate distributions is as
follows.  If the current state is $x$, we first propose a new $x^*$
according to the probabilities $S(x^*|x)$.  We then generate a series
of values $y_1,\ldots,y_{n-1}$, with $y_i$ being drawn according to
the probabilities $T_i(y_i|y_{i-1};x,x^*)$.  Let $y^*=y_{n-1}$, and
define $y_0=y$.  We accept $(x^*,y^*)$ as the new state of the Markov
chain with the following probability:
\begin{eqnarray}
  a(x,y,x^*,y^*,y_1,\ldots,y_{n-2}) 
  & \!=\! & \min\left[\, 1,\ 
  { S(x|x^*)\, \pi(x^*,y^*) \over S(x^*|x)\, \pi(x,y)}\,
  \prod_{i=1}^{n-1} {\rho_i(y_{i-1};x,x^*) \over \rho_i(y_i;x,x^*)}\,
  \right] \\[4pt]
  & \!=\! & \min\left[\, 1,\ 
  { S(x|x^*) \over S(x^*|x)}\,
    \exp\left( {1 \over n}\sum\limits_{i=0}^{n-1} E(x,y_i) \, -\,
               {1 \over n}\sum\limits_{i=0}^{n-1} E(x^*,y_i) \right)
  \right]\ \ \ \ \
\end{eqnarray}

To show that this is a valid update, I will show that the probability
in equilibrium of the chain moving from $(x,y)$ to a different
state $(x^*,y^*)$ while generating intermediate states $y_1,\ldots,y_{n-2}$
is equal to the probability of the chain moving from $(x^*,y^*)$ to
$(x,y)$ while generating intermediate states $y_{n-2},\ldots,y_1$.
Detailed balance then follows by summing over possible sequences of 
intermediate states.  The probability of moving from $(x,y)$ to
$(x^*,y^*)$ via $y_1,\ldots,y_{n-2}$ can be written\nolinebreak{} 
as\vspace*{-2pt}
\begin{eqnarray}
\lefteqn{
 \pi(x,y)\,S(x^*|x)\left[\,\prod_{i=1}^{n-1} T_i(y_i|y_{i-1};x,x^*)\,\right]
    a(x,y,x^*,y^*,y_1,\ldots,y_{n-2})}\ \ \nonumber\\[4pt]
  & \!=\! &
  \min\left[\,S(x^*|x)\,\pi(x,y) \prod_{i=1}^{n-1} T_i(y_i|y_{i-1};x,x^*),
    \right. \nonumber \\[4pt]
  && \ \ \ \ \ \ \ \left. S(x|x^*)\,\pi(x^*,y^*) \prod_{i=1}^{n-1} 
    {\rho_i(y_{i-1};x,x^*)\,T_i(y_i|y_{i-1};x,x^*) \over \rho_i(y_i;x,x^*)}\, 
  \right]\ \ \ \ \ \\[4pt]
  & \!=\! &
  \min\left[\,S(x^*|x)\,\pi(x,y) \prod_{i=1}^{n-1} T_i(y_i|y_{i-1};x,x^*),\
    S(x|x^*)\,\pi(x^*,y^*) \prod_{i=1}^{n-1} T_i(y_{i-1}|y_i;x,x^*)
  \right]\ \ \ \ \ \\[4pt]
  & \!=\! &
  \min\left[\,S(x^*|x)\,\pi(x,y) \prod_{i=1}^{n-1} T_i(y_i|y_{i-1};x,x^*),\
    S(x|x^*)\,\pi(x^*,y^*) \prod_{i=1}^{n-1} T_{n-i}(y_{i-1}|y_i;x^*,x)
  \right]\ \ \ \ \
\end{eqnarray}
If  we swap $x$ and $x^*$, $y$ and $y^*$, and 
$y_i$ and $y_{n-i-1}$, reverse the order of the
two products, and swap the arguments of $\min$, we see that this
expression is unchanged, 
showing that the reverse transition from $(x^*,y^*)$ to $(x,y)$ via
$y_{n-2},\ldots,y_1$ is equally likely.

\section{\hspace*{-7pt}Tests on simple 
                       distributions}\label{sec-tests}\vspace{-10pt}

\begin{figure}
\vspace*{-30pt}
\includegraphics{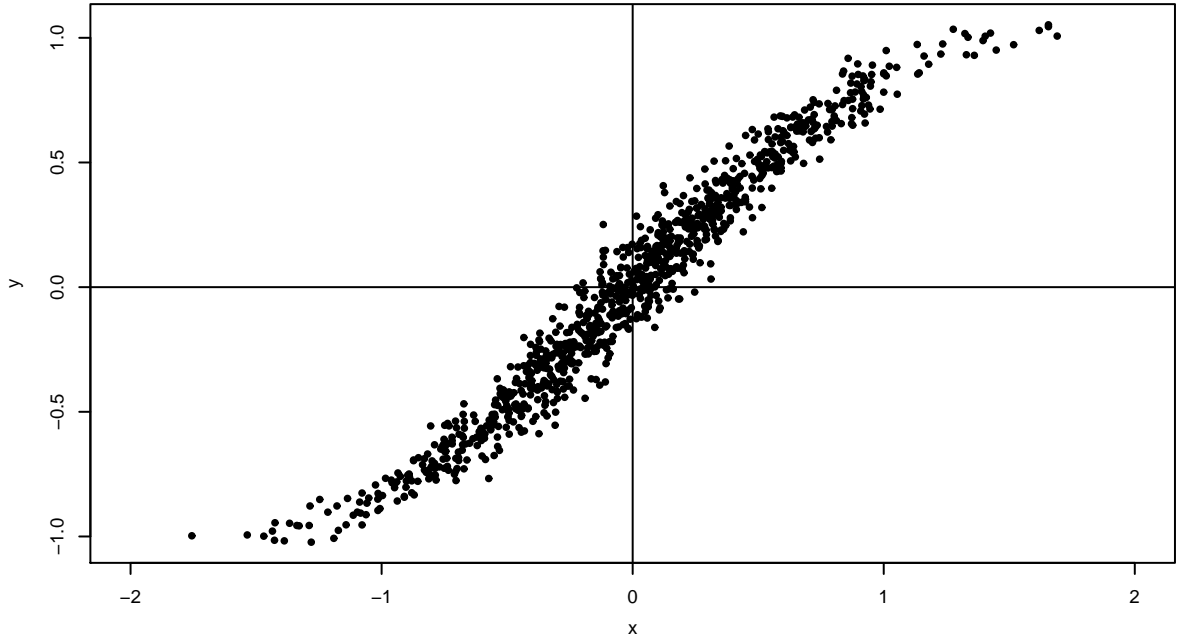}
\vspace*{-25pt}
\caption{A sample of 1000 points from the first test 
         distribution.}\label{fig-dist}
\end{figure}

I first tested the dragging method on a simple distribution in which
$x$ and $y$ are both one-dimensional, with $\pi(x,y)$ defined by the
following energy function:
\begin{eqnarray}
  E(x,y) & = & x^2\ +\ 50\,(1+x^2)^2\,(y-\sin(x))^2
\end{eqnarray}
Examination of this shows that the conditional distribution for $y$ 
given $x$ is Gaussian with mean $\sin(x)$ and standard deviation $0.1/(1+x^2)$.
From this, one can deduce that the marginal distribution for $x$ can
be obtained with an energy function of $x^2+\log(1\!+\! x^2)$.  For this
test problem, we can therefore compare performance using dragging transitions 
to the ``ideal'' performance when doing Metropolis updates based on
this marginal distribution.  Figure~\ref{fig-dist} shows a sample of
points obtained in this way, with $y$ values filled in randomly from
their conditional distribution given $x$.

For purposes of this test, we can pretend that computing $\sin(x)$ is
much slower than any of the other computations involved in evaluating
$E(x,y)$, or in the mechanics of performing Markov chain updates.
This will make $x$ a ``slow'' variable, whereas $y$ will be a ``fast''
variable.  We also pretend that we don't know that $x$ and $y$
are positively correlated.  This mimics situations in which we are 
first exploring the distribution, or in which the relationship between $x$ 
and $y$ is non-monotonic, so that no linear transformation is helpful.

\begin{figure}[p]
\hspace*{6pt}\includegraphics{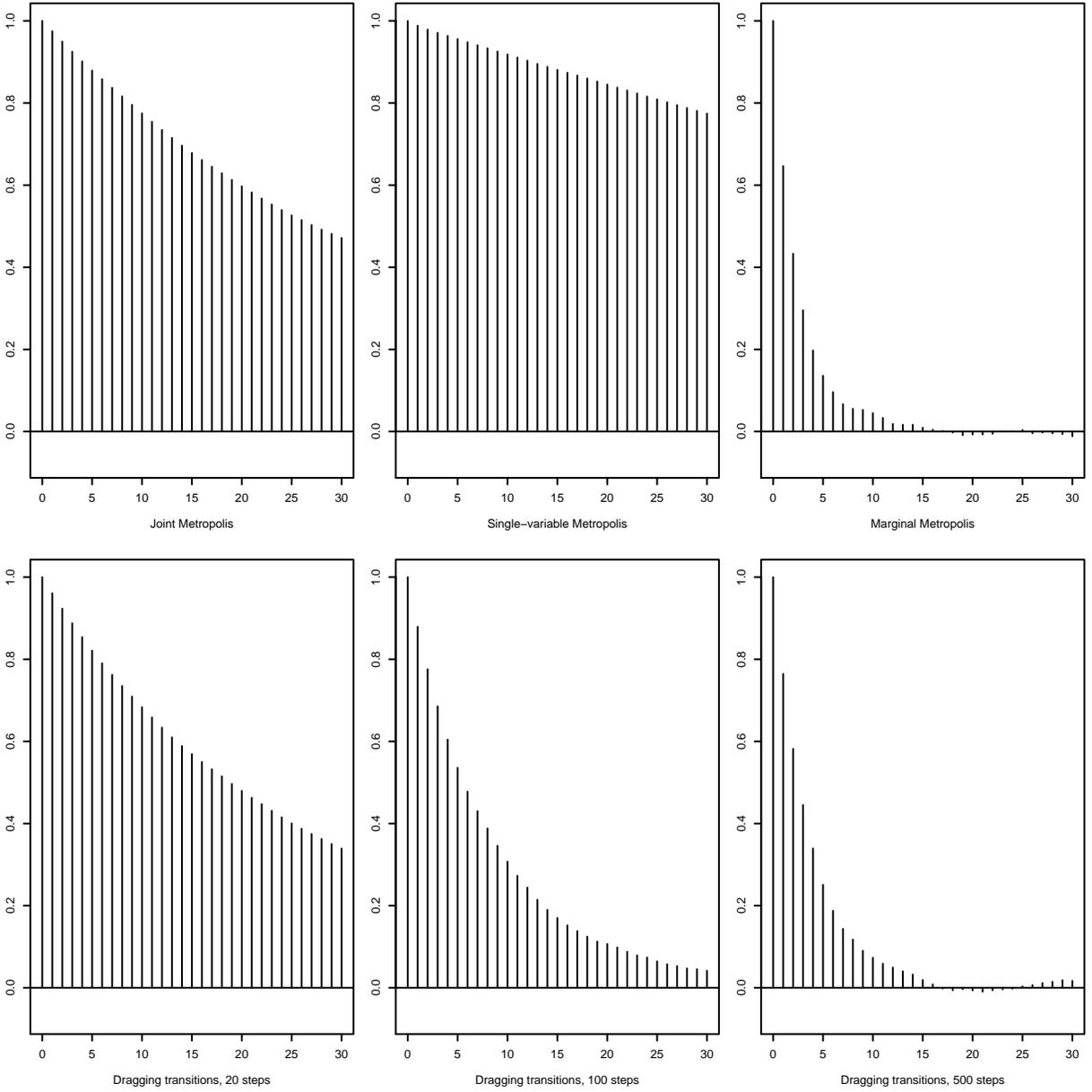}
\caption{Estimated autocorrelations for $x$ at lags up to 30 when sampling
         for the first test problem using six methods.}\label{fig-auto}
\end{figure}

Figure~\ref{fig-auto} shows the efficiency of six sampling methods
applied to this distribution, as measured by the autocorrelations for
$x$ at lags up to 30.  All the methods are based on the Metropolis
algorithm with Gaussian proposals centred on the current state.  In
all cases, the standard deviation of the Gaussian proposals was
adjusted to be approximately optimal.  All the methods require only
one slow computation of $\sin(x)$ for each iteration (for the Marginal
Metropolis method, this would be needed only when filling in $y$
values to go with the $x$ values).

In the Joint Metropolis method, the proposals change $x$ and $y$
simultaneously and independently, with the standard deviations for each
being 0.5.  The rejection rate for these proposals was 87\%.  In
the Single-variable Metropolis method, two Metropolis updates are
done each iteration, one for $x$ only, the other for $y$ only.  The
standard deviations for these proposals were both 0.25.  The rejection
rates were 59\% for $x$ and 64\% for $y$.  For the Marginal Metropolis
method, where the state consists of $x$ alone, the proposals had
standard deviation of 1.0, and the rejection rate was 47\%.  Clearly,
the Marginal Metropolis method performs much better than the other two,
though in real problems it would typically be infeasible.

The remaining plots show the autocorrelations when sampling using
updates that drag $y$ while changing $x$, with 20, 100, and 500
intermediate distributions.  For all three plots, the proposal
distribution for $x$ had standard deviation 1.0, while the proposal
distributions for $y$ during the intermediate transitions had standard
deviation 0.2.  The rejection rate for the ``inner'' updates of $y$
was around 60\% for all three runs.  The rejection rates for the
``outer'' updates of $x$ were 76\%, 63\%, and 52\% for 20, 100, and
500 intermediate distributions.  Both the rejection rate and the
autocorrelations seem to be approaching the ``ideal'' values seen
with the Marginal Metropolis method.  Provided that recomputing $E(x,y)$
after changing $y$ is around a thousand times faster than recomputing
it after changing $x$, updates for $x$ using dragging transitions will
be almost as good as updates based on the marginal distribution of $x$.

To see how sensitive these results are to the dimensionality of the 
fast parameter, I did a second test introducing another fast parameter, $z$.
The energy function used was
\begin{eqnarray}
  E(x,y) & = & x^2\ +\ 50\,(1+x^2)^2\,(y-\sin(x))^2
                  \ +\ 12.5\,(z-y)^2
\end{eqnarray}
This produces marginal distributions for $(x,y)$ and for $x$ that are the
same as for the first test.  

\begin{figure}[p]
\hspace*{6pt}\includegraphics{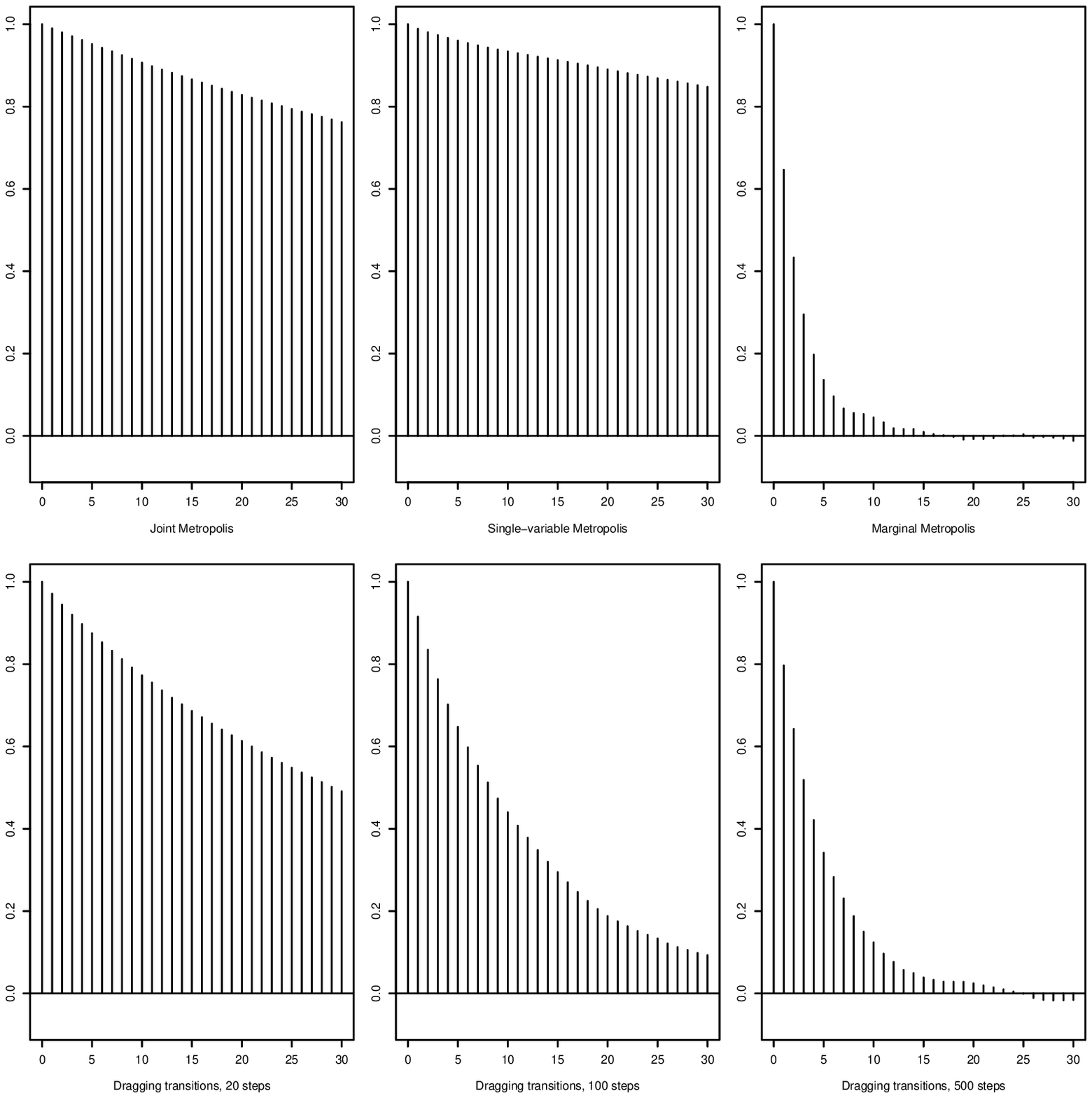}
\caption{Estimated autocorrelations for $x$ at lags up to 30 when sampling
         for the second test problem using six methods.}\label{fig-auto2}
\end{figure}

Figure~\ref{fig-auto2} shows the efficiency of the six sampling methods
applied to this distribution.  The same proposal standard deviations were
used as in the first test, except that for the Joint Metropolis updates,
the standard deviations were 0.3, producing a rejection rate of 85\%.  The
dragging transitions were done using Joint Metropolis updates for $y$ and $z$
as the inner transitions, with proposal standard deviations of 0.2.

As can be seen, all methods perform less well with the extra variable,
except for the Marginal Metropolis method, which is the same as in the
first test.  The dragging transitions are less affected, however.  The
autocorrelation time (one plus twice the sum of the autocorrelations
at all lags) when using 500 intermediate distributions increased from
approximately 7.4 to approximately 9.3 with the addition of $z$.  In
contrast, the autocorrelation time for the Joint Metropolis updates
increased from approximately 75 to approximately 205, and that for the
Single-variable Metropolis updates went from approximately 230 to
approximately 365.

The programs (written in R) used for these tests are available from my
web page.

\section*{Acknowledgements}\vspace{-10pt}

I thank Antony Lewis and Sarah Bridle for introducing me to the CMB
application, and them along with David MacKay for comments on the
manuscript.  This research was supported by the Natural Sciences and
Engineering Research Council on Canada.  I hold a Canada Research
Chair in Statistics and Machine Learning.

\section*{References}\vspace{-10pt}

\leftmargini 0.2in
\labelsep 0in

\begin{description}
\itemsep 2pt

\item
Diggle, P.~J., Tawn, J.~A., and Moyeed, R.~A.\ (1998) ``Model-based
  geostatistics'', {\em Applied Statistics}, vol.~47, pp.~299-350.

\item
  Hastings, W.~K.\ (1970) ``Monte Carlo sampling methods using Markov chains 
  and their applications'', {\em Biometrika}, vol.~57, pp.~97-109.

\item
  Lewis, A.\ and Bridle, S.\ (2002) ``Cosmological parameters from CMB
    and other data: a Monte-Carlo approach'', 
    \texttt{http://arxiv.org/abs/astro-ph/0205436}.

\item
  Metropolis, N., Rosenbluth, A.~W., Rosenbluth, M.~N., Teller, A.~H., 
  and Teller, E.\ (1953) ``Equation of state calculations by fast computing 
  machines'', {\em Journal of Chemical Physics}, vol.~21, pp.~1087-1092.

\item
  Neal, R.~M.\ (1996) ``Sampling from multimodal distributions using
  tempered transitions'', {\em Statistics and Computing}, vol.~6, pp.~353-366.

\item
  Neal, R.~M.\ (1999) ``Regression and classification using Gaussian process
      priors'' (with discussion), in J.~M.~Bernardo, {\em et al} 
      (editors) {\em Bayesian Statistics 6}, Oxford University Press, 
      pp.~475-501.

\end{description}

\end{document}